\input amstex
\documentstyle{amsppt}
\magnification=\magstep1

\TagsOnRight \NoBlackBoxes

\hoffset1 true pc \voffset2 true pc

\hsize36 true pc \vsize52 true pc

\tolerance=2000
\def\m1{^{-1}}
\def\ov1{\overline}
\def\gp#1{\langle#1\rangle}
\def\cry#1{\operatorname{\frak{Crys}}(#1)}

\catcode`\@=11
\def\logo@{}
\catcode`\@=\active

\topmatter

\title
Extensions of the representation modules of a prime order group
\endtitle
\author
V.~A.~Bovdi,  V.~P.~Rudko
\endauthor
\leftheadtext\nofrills{ V.~A.~Bovdi,  V.~P.~Rudko} \rightheadtext
\nofrills {Extensions of the representations modules }
\dedicatory
Dedicated to Professor {\it K\'alm\'an Gy\H ory}  on his 65th
birthday
\enddedicatory
\abstract For the ring $R$ of integers of a ramified  extension of
the field of $p$-adic numbers and a cyclic group $G$ of prime
order $p$ we study the extensions of the additive groups of
$R$-representations modules of $G$ by the group $G$.
\endabstract
\thanks
The research was supported by OTKA  No.T 037202 and  No.T 038059
\endthanks
\address
\hskip-\parindent {\rm  V.A.~Bovdi\newline Institute of
Mathematics, University of Debrecen\newline P.O.  Box 12, H-4010
Debrecen, Hungary\newline Institute of Mathematics and
Informatics, College of Ny\'\i regyh\'aza\newline S\'ost\'oi \'ut
31/b, H-4410 Ny\'\i regyh\'aza, Hungary}
\newline
E-mail: vbovdi\@math.klte.hu
\bigskip
\hskip-\parindent {\sl V.P.~Rudko}
\newline
Department of Algebra,  University of Uzhgorod
\newline
88 000, Uzhgorod, Ukraine
\newline
E-mail: math1\@univ.uzhgorod.ua
\endaddress
\endtopmatter

\document
Let $\frak{T}$ be the field of fractions of a principal ideal
domain $R$, $F$ a field which contains $R$, let $G$ be a finite
group and $\Gamma$ a matrix $R$-representation of $G$. Let $M$ be
an $RG$-module, which affords  the $R$-representation $\Gamma$ of
$G$, and $FM=F\otimes_RM$ the smallest linear space over $F$ which
contains $M$ and $\widehat{M}=FM^+/M$, the factor group of the
additive group of the space $FM$ by the additive group of $M$.
Clearly, the group $\widehat{M}$ and the space $FM$ are
$RG$-modules. Put $\widehat{F}=F^+/R$.

Let $\frak{f}:G\to\widehat{M}$ be a $1$-cocycle of $G$ with value
in  $\widehat{M}$, i.e.
$$
\frak{f}(xy)=x\frak{f}(y)+\frak{f}(x), \quad\quad\quad\quad
(x,y\in G).
$$
Define $[g,x]$  by $\left(\smallmatrix g& x\\
0&1\\\endsmallmatrix\right)$ and set
$$
\cry{G,M,\frak{f}}=\{\; [g,x]\;  \mid \; g\in G, \quad x\in
\frak{f}(g)\; \},
$$
where $x$ runs over the cosets $\frak{f}(g)\in \widehat{M}$ for
any $g\in G$.

Clearly, $\cry{G,M,\frak{f}}$ is a group, where the multiplication
is the usual matrix multiplication.  Of course $ K_1=\{ [e,x] \mid
e \; \text{is the unit element of} \; G, \quad x\in \frak{f}(e)
\}$ is a normal subgroup of $\cry{G,M,\frak{f}}$ such that
$K_1\cong M^+$ and $\cry{G,M,\frak{f}}/K_1\cong G$. The group
$\cry{G,M,\frak{f}}$ is an extension of the additive group of the
$RG$-module $M$ by $G$.

We are using the terminology of the theory of group
representations \cite{1}.

A $1$-cocycle $\frak{f}: G\to \widehat{M}$ is called {\it
coboundary}, if  there exists an $x\in FM$ such  that
$\frak{f}(g)=(g-1)x+M$ for every $ g\in G$. The $1$-cocycles
$\frak{f}_1: G\to \widehat{M}$ and $\frak{f}_2: G\to \widehat{M}$
are called {\it cohomologous} if $\frak{f}_1-\frak{f}_2$ is a
coboundary. Let $H^1(G,\widehat{M})$ be the first cohomology
group. Clearly, each element of $H^1(G,\widehat{M})$ defines a
class of equivalence of groups.

If the $1$-cocycles $\frak{f}_1,\frak{f}_2$ are cohomologous, then
$\cry{G,M,\frak{f}_1}$ and $\cry{G,M,\frak{f}_2}$ are isomorphic.
This isomorphism is called {\it equivalence} and these groups  are
called equivalent. In particular, the group $\cry{G,M,\frak{f}}$
is split (i.e. $\cry{G,M,\frak{f}}=M\rtimes G$) if and only if
$\frak{f}$ is coboundary.

The {\it dimension} of the group $\cry{G,M,\frak{f}}$ is called
the $R$-rank of the $R$-module $M$. (Note that $M$ is a free
$R$-module of finite rank.)   The group $\cry{G,M,\frak{f}}$ is
called {\it irreducible} ({\it indecomposable}), if $M$ is an
irreducible (indecomposable) $RG$-module and the $1$-cocycle
$\frak{f}$ is not cohomologous  to  zero.

The group $\cry{G,M,\frak{f}}$ is  {\it non-split}, if the
$1$-cocycle $\frak{f}$ defines a nonzero element of
$H^1(G,\widehat{M})$.

Note that the properties of the group $\cry{G,M,\frak{f}}$ were
studied in \cite{5,6,8}, in the cases when $R$ is either the ring
of rational integers $\Bbb Z$, or  the $p$-adic integers $\Bbb
Z_{p}$, or the localization  $\Bbb Z_{(p)}$  of $\Bbb Z$ at $p$.

Let $G=\gp{a\mid a^p=1}$ be the cyclic group of prime order $p$,
$R$ the ring of integers of the ramified  finite extension  $T$ of
the field of $p$-adic numbers. We calculate the group
$H^1(G,\widehat{M})$ for some module $M$ of an  indecomposable
$R$-representation of $G$.

Let $\Phi_p(x)=x^{p-1}+\cdots+x+1$ be a cyclotomic polynomial of
degree $p$ and let $\eta(x)$ be a divisor of $\Phi_p(x)$ over the
field $\frak{T}$ with $deg(\eta(x))<p-1$ (provided that   such
nontrivial polynomial exists).

\proclaim {Lemma 1} Let $M_1$ and $M_2$ be $RG$-modules which
afford  an $R$-rep\-re\-sen\-tation $\Gamma$ of $G=\gp{a\mid
a^p=1}$.
\itemitem{(i)} If $M_1\cong M_2$ then
$H^1(G,\widehat{M_1})\cong H^1(G,\widehat{M_2})$.
\itemitem{(ii)} If  the matrix $\Gamma(a)$ does not have $1$ as eigenvalue, then $H^1(G,\widehat{M_1})$ is trivial.
\endproclaim
\demo{Proof} See \cite{1}. \qquad\qed \enddemo

\proclaim {Theorem  1} Let $G=\gp{a\mid a^p=1}$ and
$M_\eta=\eta(a)RG$. Then the $RG$-module $M_\eta$ is
indecomposable and
$$
H^1(G,\widehat{M_\eta})\cong R/(\eta(1)R),
$$
where $R/(\eta(1)R)$ is the additive group of the factor ring of
$R$ by the ideal  $\eta(1)R$.
\endproclaim
\demo{Proof} Let $t\in R$ be a prime element  and
$\ov1{R}=R/(tR)$. Then in  $\ov1{R}$ we have that
$$
x^p-1=(x-1)^p;\qquad \eta(x)=(x-1)^n;\qquad
(x^p-1)\eta\m1(x)=(x-1)^{p-n},\tag1
$$
where $n=deg(\eta(x))$.

Put\quad   $\eta_1(x)=(x^p-1)\eta\m1(x)$.\quad  Then $M_\eta$ and
$RG/(\eta_1(a)RG)$ are isomorphic as $RG$-modules. If
$\ov1{M_\eta}=M_\eta/(tM_\eta)$, then by (1) follows that
$\ov1{M_\eta}$ is a root subspace  of the linear operator $a$ over
$\ov1{R}$. It is easy to see that $\ov1{M_\eta}$ is not
decomposable into a direct sum of  invariant subspaces. It follows
that $M_\eta$ is an indecomposable $RG$-module. Clearly \quad
$FM_\eta=\eta(a)FG=\eta(a)F+(a-1)FM_\eta$\quad  and the group
$\widehat{M_\eta}=FM_\eta^+/M_\eta$ is isomorphic to  a direct sum
of groups\quad  $\eta(a)(F^+/R)+(a-1)\widehat{M_\eta}$. This means
that in the class  of $1$-cocycles there is a cocycle $\frak{f}:
G\to \widehat{M_\eta}$ such that
$$
\frak{f}(a)=\lambda \eta(a)+M_\eta, \quad\quad\quad (\lambda\in
F).
$$
Moreover, from $\frak{f}(a^p)=0$ (in $\widehat{M_\eta}$) it
follows that if $\omega=a^{p-1}+\cdots+a+1$, then
$$
\omega\cdot\frak{f}(a)=\lambda\cdot\eta(1)\omega\in M_\eta
$$
if and only if $\lambda\eta(1)\in R$. Therefore,  $H^1(G,
\widehat{M_\eta})$ is isomorphic to the subgroup \quad $\{\;
\lambda+R\; \mid \;\lambda\in F,\;\lambda\cdot\eta(1)\in R\;\}$ of
$F/R$ and
$$
\{\; \lambda+R\in F/R\; \mid \; \lambda\cdot\eta(1)\in R\;\}\cong
R/(\eta(1)R).
$$
\rightline{\text{\qed}}
\enddemo

 \proclaim {Corollary  1} Let $\alpha\in R$ and
$\eta(1)R=t^sR$, where $t$ is a prime element of $R$. Put
$$
{K}_\alpha(G,M_\eta)=\bigg\langle{\;\left(\smallmatrix e& \quad m\\
0&\quad 1\\\endsmallmatrix\right),\quad
\left(\smallmatrix a&\quad \alpha t^{-s}\eta(a)\\
0&1\\\endsmallmatrix\right)\quad \mid \quad m\in
M_\eta}\;\bigg\rangle,
$$
where $\alpha$ runs over the representative elements of the cosets
of  $R/(t^sR)$. Up to equivalence,  the groups
$K_\alpha(G,M_\eta)$ give all extensions of the additive group of
the $RG$-module $M_\eta$ by the group $G$.
\endproclaim

\bigskip

Suppose $p=t^d\theta$, where $d>1$ is the ramification index  and
$\theta$ is a unit in $R$. Set
$$
{\frak X}_{ji}=t^jRG+(a-1)^iRG,\quad\quad \quad  (1\leq j<d,\;
1\leq i<p).
$$

\proclaim {Theorem 2} The module ${\frak X}_{ji}$ is an
$RG$-module affording an indecomposable $R$-representation of $G$
and
$$
H^1(G,\widehat{{\frak X}_{ji}})\cong R/(t^{d-j}R).
$$

\endproclaim
\demo{Proof} Suppose that the $RG$-module ${\frak X}_{ji}$ is
decomposable into a direct sum of $RG$-submodules. Then
$t^j=u_1+u_2$, where $u_1,u_2$ are nonzero elements of $RG$ with
$u_1u_2=0$. Thus $e_1=t^{-j}u_1$ is an idempotent. Since the trace
$tr(e_1)$ of $e_1$ is a rational number (see Theorem 3.5, \cite{4}
p.21) of the form $rp\m1$ \quad ($1\leq r\leq p$), we get
$t^jrp\m1\in R$, which is impossible for $j<d$. This contradiction
proves the indecomposability   of ${\frak X}_{ji}$. Clearly \quad
$F{\frak X}_{ji}=FG=F+(a-1)FG$. \quad  Therefore, in each class of
$1$-cocycles there is  a cocycle $\frak{f}: G \longrightarrow
\widehat{{\frak X}}_{ji}$ such that $ \frak{f}(a)=\lambda+{\frak
X}_{ji}$, where $\lambda\in F$ with $\lambda\omega\in {\frak
X}_{ji}$. It follows that $\lambda p=t^j\alpha$, where $\alpha\in
R$ and
$$
H^1(G,\widehat{{\frak X}}_{ji})\cong \{\lambda+R\;\mid\;
\lambda\in F, \;\lambda t^{d-j}\in R\}
$$
is a subgroup of $F/R$.

\rightline{\text{\qed}}
\enddemo

Set
$$
K_\alpha(G,{\frak X}_{ji})=\bigg\langle{\;\left(\smallmatrix e& \quad m\\
0&\quad 1\\\endsmallmatrix\right),\quad
\left(\smallmatrix a&\quad \alpha t^{j-d}\\
0&1\\\endsmallmatrix\right)\quad \mid \quad m\in {\frak
X}_{ji}}\;\bigg\rangle,
$$
where $\alpha$ runs over the representative elements of the cosets
of  $R/(t^{d-j}R)$.

\proclaim {Corollary  2} The groups $K_\alpha(G,{\frak X}_{ji})$
 give all extensions of the additive group of the
$RG$-module ${\frak X}_{ji}$ by $G$.
\endproclaim

\bigskip

\proclaim {Lemma 2} The set $\textstyle \{\; {\frak X}_{ji}\; \mid
j=1,\ldots,d-1;\; i=1,\ldots,\frac{p-1}{2}\}$ consists of pairwise
non-isomorphic modules.
\endproclaim
\demo{Proof} Let us consider an  indecomposable $\ov1{R}G$-module
$V_i=\ov1{R}G\big/\big((a-1)^i\ov1{R}G\big)$, where
$\ov1{R}=R/(tR)$ and $1\leq i\leq p$. It is easy to check that the
elements
$$
u_1=t^j,\ldots, u_i=t^j(a-1)^{i-1},\quad u_{i+1}=(a-1)^{i},
\ldots, u_p=(a-1)^{p-1}\tag2
$$
form an $R$-basis in ${\frak X}_{ji}$ and
$$
\Phi_p(x)-(x-1)^{p-1}=p\theta(x),\tag3
$$
where $\theta(x)\in \Bbb Z[x]$, $deg(\theta(x))\leq p-2$. Note
that since $\theta(1)=1$, it  follows that $\theta(a)$ is a unit
in the group ring $RG$. Using the identity
$$
xy-1=(x-1)(y-1)+(x-1)+(y-1),
$$
from (3) we obtain that
$$
(a-1)^p=p(a-1)\cdot\big(\alpha_0+\alpha_1(a-1)+\cdots+\alpha_{p-2}(a-1)^{p-2}\big),\tag4
$$
where $\alpha_0,\alpha_1,\ldots,\alpha_{p-2}\in \Bbb Z$. Since
$p=t^d\theta=t(t^{d-1}\theta)$, from (4) we get
$$
(a-1)^p=(a-1)u_p=tm,\tag5
$$
where $m\in {\frak X}_{ji}$. According to (2) \quad
$(a-1)u_i=t^ju_{i+1}$,  and from (5) we obtain that the
$RG$-module $\ov1{\frak X}_{ji}={\frak X}_{ji}/(t{\frak X}_{ji})$
is isomorphic to a direct sum $V_i\oplus V_{p-i}$ of
indecomposable $\ov1{R}G$-modules, so by Theorem 2 and Lemma 1 the
proof is complete.

\rightline{\text{\qed}}
\enddemo

Let $n>1$ be the degree  of a  divisor of $\Phi_p(x)$, which is
irreducible over $R$. We consider the following $RG$-modules:
$$
{\frak U}_{ji}=t^j(a-1)RG+(a-1)^{s+1}RG,\quad\quad\quad (1\leq
j<d,\quad  1\leq s<n).
$$
It is easy  to check that the $RG$-module ${\frak U}_{ji}$
satisfies the condition (ii) of Lemma 1, so $H^1(G,\widehat{{\frak
U}_{ji}})=0$.

Let ${\frak Z}_{js}$ be a submodule of the free module \quad
$RG^{(2)}=\{(x,y)\mid x,y\in RG\}$\quad  of rank $2$, which
consists of the  solutions $(x,y)$ of the  equality
$$
t^j(a-1)x+(a-1)^{s+1}y=0.\tag6
$$

\proclaim {Lemma 3} Let $\omega=\Phi_p(a)$ and set $u_1=[0,
\omega]$, $u_2=[(a-1)^s, -t^{j}]$ and $u_3=[
t^{-j}(\omega-(a-1)^{p-1}), (a-1)^{p-s-1}]$. Then ${\frak Z}_{js}$
is an $RG$-module generated by $u_1,u_2,u_3$.
\endproclaim

\demo{Proof} Clearly, $u_1,u_2,u_3\in {\frak Z}_{js}$. Let
$u=[x,y]$ be an arbitrary element of ${\frak Z}_{js}$. If $x=0$
then $u=u_1$. Suppose $x\not=0$. By  substraction  of the elements
of $RGu_3$ from $u$  we  obtain that
$y=\gamma_0+\gamma_1(a-1)+\cdots+\gamma_{p-s-2}(a-1)^{p-s-2}$\quad
($\gamma_r\in R$).  By (6)
$$
t^{j}(a-1)x+\big(\gamma_0+\gamma_1(a-1)+\cdots+\gamma_{p-s-2}(a-1)^{p-s-2}\big)\cdot(a-1)^{s+1}=0,
$$
which is possible if and only if\quad
$\gamma_0\equiv\cdots\equiv\gamma_{p-s-2}\equiv 0\pmod{t^{j}}$.
\quad  Now, since  $u$ is an element of $RGu_2$, we  obtain that
$y=0$. Then\quad  $t^j(a-1)x=0$\quad  which implies
$x=\alpha\omega$ \; ($\alpha\in R$) and
$u=\alpha(t^ju_3-(a-1)^{p-s-1}u_2)$.

\rightline{\text{\qed}}\enddemo

\proclaim {Theorem  3} The $RG$-module ${\frak Z}_{js}$ is
indecomposable. Moreover,
$$
H^1(G,\widehat{{\frak Z}}_{js})\cong R/(t^{d}R)\oplus R/(t^{d-j}R)
$$
and the $RG$-modules ${\frak X}_{js}$ are pairwise non-isomorphic.
\endproclaim
\demo{Proof} It is easy to see that
$$
u_1=t^j(a-1),\ldots, u_{i-1}=t^j(a-1)^{i-1},\quad u_{i}=(a-1)^{i},
\ldots, u_{p-1}=(a-1)^{p-1}
$$
form an $R$-basis in the $RG$-module ${\frak U}_{js}$ and
$$
\ov1{{\frak U}}_{js}={\frak U}_{js}/(t{\frak U}_{js})\cong
V_s\oplus V_{p-s-1}.
$$
Since $s<n$,  it follows that the  $RG$-module ${\frak U}_{js}$ is
indecomposable. Moreover, it follows that the $RG$-modules ${\frak
U}_{js}$ are pairwise  non-isomorphic and $RG$-modules ${\frak
Z}_{js}$, ${\frak U}_{js}$  and  $ RG^2$ form an exact sequence
$$
0 @>{ }>>\frak Z_{js}@>{}>>  RG^{(2)} @>{}>>\frak U_{js}@>{}>>0.
$$
Therefore, $\frak Z_{js}$ is the kernel of a minimal projective
covering of the indecomposable $RG$-module $\frak U_{js}$, so
$\frak Z_{js}$ is also indecomposable.

\proclaim {Lemma 4} Let $\widehat{\frak Z}_{js}=(F\frak
Z_{js})^+/\frak Z_{js}$, $\widehat{F}=F^+/R$ and
$M=(a-1)\widehat{\frak Z}_{js}$. Then
$$
\widehat{\frak Z}_{js}/M=\widehat{F}\nu_1+\widehat{F}\nu_2,
$$
where $\nu_1=[0,\omega]+M$,\quad  $\nu_2=[\omega,0]+M$\quad
and\quad  $a\nu_1=\nu_1$,\quad  $a\nu_2=\nu_2$.
\endproclaim
\demo{Proof} Clearly, $ax=x$\quad  ($x\in \widehat{\frak
Z}_{js}/M$) and $\widehat{F}\nu_1=\widehat{F}[0,\omega]+M\in
\widehat{\frak Z}_{js}/M$. Moreover,
$$
\split
\omega\widehat{F}u_3+M=\widehat{F}[t\m1p\omega,0]+M=&\widehat{F}(tp\m1)[t\m1p\omega,0]+M\\
=&\widehat{F}[\omega,0]+M.
\endsplit
$$
By  analogy
$$
\omega\widehat{F}u_2+M=\widehat{F}[0,-t^j\omega]+M=
\widehat{F}[0,\omega]=\widehat{F}\nu_1.
$$
Therefore $\widehat{\frak
Z}_{js}/M=\widehat{F}\nu_1+\widehat{F}\nu_2$. {\text{\qed}}
\enddemo
\bigskip

From Lemma 4 it follows that  each class of $1$-cocycles of the
group $G$ with  values in the group $\widehat{\frak
Z}_{js}=(F{\frak Z}_{js})^+/{\frak Z}_{js}$  contains a
$1$-cocycle $\frak{f}$ such that
$$
\frak{f}(a)=\alpha[0,\omega]+\beta[\omega,0]+{\frak Z}_{js},
$$
where $\alpha,\beta\in F$ and \quad
$\omega\big(\alpha[0,\omega]+\beta[\omega,0]\big)\in {\frak
Z}_{js}$.  This condition holds if and only if
$\alpha{p},\beta{p}\in R$. Moreover
$$
\alpha[0,\omega]+\beta[\omega,0]\in {\frak
Z}_{js}+(a-1)\widehat{\frak Z}_{js}
$$
if and only if $\alpha\in R$ and $\beta\in t^{-j}R$. Using
properties of the $1$-cocycle $\frak{f}$ it is easy to show that
the two $1$-cocycles $\frak{f}_j$  ($j=1,2$):
$$
\frak{f}_1(a)=\alpha_1[0,\omega]+\beta_1[\omega,0]+{\frak
Z}_{js},\qquad
\frak{f}_2(a)=\alpha_2[0,\omega]+\beta_2[\omega,0]+{\frak Z}_{js}
$$
are cohomologous if and only if
$$
p\alpha_1\equiv p\alpha_2\pmod{t^d} \qquad\quad
\text{and}\qquad\quad p\beta_1\equiv p\beta_2\pmod{t^{d-j}},
$$
where  $\alpha_j,\beta_j\in F$,\quad  $p\alpha_j, p\beta_j\in R$.
Note that  $p=t^d\theta$.

It follows that the map\quad  $\frak{f}\longmapsto
\big(p\alpha+t^dR,\; p\beta+t^{d-j}R\big)$\quad  gives the
isomorphism
$$
H^1(G,\widehat{{\frak Z}}_{js})\cong R/(t^{d}R)\oplus
R/(t^{d-j}R).
$$
Therefore, according to (ii) of Lemma 1, the $RG$-modules
$\frak{Z}_{js}$ \quad ($1\leq j<d$) are pairwise non-isomorphic.

\rightline {\text{\qed}}
\enddemo
Now, using the description of $1$-cocycles it is easy to prove the
following
\proclaim {Corollary  3} Put
$$
K_{\alpha,\beta}(G,\frak{Z}_{js})=\bigg\langle{\;\left(\smallmatrix e& \quad m\\
0&\quad 1\\\endsmallmatrix\right),\quad
\left(\smallmatrix a&\quad \alpha t^{-d}[0,\omega]+\beta t^{-d}[\omega,0]\\
0&1\\\endsmallmatrix\right)\quad \mid \quad m\in
Z_{js}}\;\bigg\rangle,
$$
where $\alpha$ and $\beta$ independently run over the
representative elements of the  cosets $R/(t^{d}R)$ and
$R/(t^{d-j}R)$, respectively. Up to equivalence,  the groups
$K_{\alpha,\beta}(G,\frak{Z}_{js})$ give all extensions of the
additive group of the $RG$-module $\frak{Z}_{js}$ by the group
$G$.
\endproclaim

If $R$ is the quadratic extension of the ring  of $p$-adic
integers, then the $R$-representations of $G$ were described by
P.M.~Gudivok (see \cite{7}).

\bigskip
Finally, we have the following result
 \proclaim {Theorem 4} Let
$\Phi_p(x)$ be decomposable into the product of at least two
irreducible polynomials over $R$. Then the dimensions of the
non-split indecomposable groups $\cry{G,M,\frak{f}}$ are
unbounded.
\endproclaim

\demo{Proof} Let $\Phi_p(x)=\eta_1(x)\cdots \eta_k(x)$\quad
($k>2$) be a decomposition into a product of polynomials
irreducible over $R$ and suppose that
$$
\eta_1(x)=x^n-\alpha_{n-1}x^{n-1}-\cdots-\alpha_{1}x-\alpha_{0}\in
R[x].
$$
Note that $deg(\eta_1(x))=deg(\eta_2(x))=\cdots=deg(\eta_k(x))=n$
\quad  and\quad   $kn=p-1$.
\bigskip

We will use the technic of the integral representation of finite
groups, which was developed by S.D.~Berman and P.M.~Gudivok  in
\cite{2, 3, 7}.

Let $\varepsilon$ be a primitive $p$th root of unity such that
$\eta_1(\varepsilon)=0$ and let $r_j$ be a natural number, such
that $\varepsilon_j=\varepsilon^{r_j}$  is a root of the
polynomial $\eta_j(x)$, where  $ r_1=1$  and $j=1,\ldots,k$. Let $
\widetilde{\varepsilon}= \left(\smallmatrix
0&\cdots&0&\alpha_0\\
1&\cdots&0&\alpha_1\\
\vdots&\ddots&\vdots&\vdots\\
0&\cdots&1&\alpha_{n-1}\\
\endsmallmatrix\right)
$ be the comparing matrix of  $\eta_1(x)$.

The following $R$-representations  of $G=\gp{a\mid a^p=1}$ are
irreducible:
$$
 \delta_0:\;  a\mapsto 1;\qquad \delta_1: \; a\mapsto
\widetilde{\varepsilon};\qquad \delta_j: \; a\mapsto
\widetilde{\varepsilon_j}=
\widetilde{\varepsilon_j}^{r_j},\qquad\qquad (j=2,\ldots,k).
$$
Note that the module which affords  representation $\delta_1$ is
$R[\varepsilon]$ with $R$-basis $1,\varepsilon,\ldots,
\varepsilon^{n-1}$.
\enddemo

Let  $m\in\Bbb N$. Define  the following $R$-representation of
$G=\gp{a}$ of  degree $(3n+1)m$:
$$
\Gamma_{m}:\quad a\mapsto  \quad\left(\smallmatrix \Delta_{1m}(a)\quad& U_m(a)\\
0\quad&\Delta_{2m}(a)\\\endsmallmatrix\right),
$$
where
\item{} $ \Delta_{1m}(a)=
\quad\delta_0^{(m)}(a)+\delta_1^{(m)}(a)=\quad
\left(\smallmatrix E_m\otimes \delta_0(a)\quad& 0\\
0\quad&E_m\otimes \delta_1(a)\\
\endsmallmatrix\right)$;
\smallskip
\item{} $\Delta_{2m}(a)=
\quad\delta_2^{(m)}(a)+\delta_3^{(m)}(a)=\quad
\left(\smallmatrix E_m\otimes \delta_2(a)\quad& 0\\
0\quad&E_m\otimes \delta_3(a)\\
\endsmallmatrix\right)$;
\smallskip
\item{}$U_{m}(a)= \quad\left(\smallmatrix E_m\otimes u & \quad J_m(1)\otimes u\\
E_m\otimes \overline{u}&\quad E_m\otimes \overline{u}\\
\endsmallmatrix\right)$;
\smallskip
\item{} $u=(0,0,\ldots,0,1)$\quad   defines a nonzero element of
$Ext(\delta_0,\delta_j)$;
\smallskip
\item{}$J_m(\lambda)$ \quad is a Jordan block of degree $m$ with $\lambda$
in the main diagonal;
\smallskip
\item{}$\overline{u}$\; is a matrix in which  the first row is
$(0,\ldots,0,1)$ and all other rows are zero. The matrix
$\overline{u}$ defines a nonzero element of the group
$Ext(\delta_1,\delta_j)$, where $j=2,3$;
\smallskip
\item{} $E_m$ is the unity matrix of degree $m$.

\proclaim {Lemma 5}(see \cite{2,3})  $\Gamma_m$ is an
indecomposable $R$-representation of $G$.
\endproclaim

Let $\frak{W}_m=R^l$ be a module of $l$-dimension vectors over $R$
affording the  $R$-representation $\Gamma_m$. Put
$\widehat{F}=F^+/R$,
$\widehat{\frak{W}}_m=F\frak{W}_m^+/\frak{W}_m$. Clearly
$\widehat{F}^l\cong \widehat{\frak{W}}_m$. Define $\tau:
F\longrightarrow F^n$  by
$$
\tau(w)=w\big(\alpha_0,\quad \alpha_0+\alpha_1,\quad
\alpha_0+\alpha_1+\alpha_2,\quad  \ldots, \quad
\alpha_0+\cdots+\alpha_{n-2},\quad  1\big),\tag7
$$
where the $\alpha_j$ are coefficients of $\eta_1(x)$ and $w\in F$.

\proclaim {Lemma 6} (i) Each $1$-cocycle  of $G=\gp{a\mid a^p=1}$
at $\widehat{\frak{W}}_m$ is cohomologous to a cocycle $\frak{f}$,
such that
$$
\frak{f}(a)=(X,0,\ldots,0)+\frak{W}_m,
$$
where $ X\in F^m$ and $pX=0$ in $\widehat{F}^m$ (i.e. $pX\in
R^m$).
\itemitem{(ii)} Let $z\in F^n$ such that
$(\widetilde{\varepsilon}-E_n)z=0$ in $\widehat{F}^n$. Then
$z=\tau(w)\pmod{R^n}$, with  $w\in F$ such that \; $\eta_1(1)w=0$
in $\widehat{F}$.
\itemitem{(iii)} If $V=R/(\frac{p}{\eta(1)}R)$ is the residual  of ring $R$ by
the ideal $(\frac{p}{\eta(1)}R)$, then
$H^1(G,\widehat{\frak{W}}_n)\cong V^m$.

\endproclaim
\demo{Proof} (i)\;  follows from (ii) of Lemma 1. (ii)\; is easy
to check.
\newline
(iii) \; by (i) we can put  $\frak{f}(a)=(X,0,0,0)$ and
$\frak{g}(a)=(Y,0,0,0)$, where $X=(x_1,\ldots,x_n)$,
$Y=(y_1,\ldots,y_n)$ and  $pX=pY=0$. Note that all the equalities
considered here are understood modulo the group $R$. Suppose that
these $1$-cocycles are cohomologous and $Z\in F^l$ is  such that
$$
\big(\Gamma_m(a)-E_l\big)Z+\frak{f}(a)=\frak{g}(a).\tag8
$$
Put $Z=(Z_1,Z_2,Z_3,Z_4)$, where $Z_1\in F^m$ and $Z_2,Z_3,Z_4$
are $m$-dimensional vectors, with  $i$-components belong to $F^n$,
and  denoted  by $Z^i_2,Z^i_3$ and $Z^i_4$, respectively. By (8)
we get
$$
(E_m\otimes u)Z_3+(J_m\otimes u)Z_4+X=Y,\tag9
$$$$
\big(E_m\otimes(\widetilde{\varepsilon} -E_n)\big)Z_2+(E_m\otimes
\overline{u})(Z_3+Z_4)=0,\tag10
$$$$
\big(E_m\otimes(\widetilde{\varepsilon}_2 -E_n)\big)Z_3=0,\qquad
\big(E_m\otimes(\widetilde{\varepsilon}_3 -E_n)\big)Z_4=0.\tag11
$$
From (11) and by (ii) we have
$$
Z_3=\big(\tau(v_1),\ldots,\tau(v_m)\big),\qquad
Z_4=\big(\tau(u_1),\ldots,\tau(u_m)\big),\tag12
$$
where $u_j, v_j\in F$,\quad   $\tau$  is from (7) and
$$
\eta_1(1)u_j=\eta_1(1)v_j=0.\tag13
$$
Clearly, the equality (10) consists of  $m$ matrix equalities of
the form
$$
(\widetilde{\varepsilon}-E_n)Z_2^i+\overline{u}\tau(w)=0,\tag14
$$
where  $Z_2^i\in F^n$ is the $i$th component of $Z_2$,\quad
$i=1,\ldots,m$ and $w\in F$. Since  $u\tau(w)=w$ and
$\overline{u}\tau(w)=(w,0,\ldots,0)$,  when all the rows of (14)
are added together we obtain
$$
-\eta(1)Z_2^n+w=0, \tag15
$$
where $Z_2^n$ is the last component of the vector $Z_2$. According
to (12) and (15), (10) gives the equalities
$$
-\eta(1)z_j+v_j+u_j=0,\qquad (j=1,\ldots,m)\tag16
$$
where $z_j$ are some components of $Z_2$. From (9)
$$
\split
v_j+u_j+u_{j+1}+x_j&=y_j,\qquad\quad(j=1,\ldots,m-1)\\
v_m+u_m+x_m&=y_m,\\
\endsplit\tag17
$$
where $X=(x_1,\ldots,x_n)$, $Y=(y_1,\ldots,y_n)$ and  $pX=pY=0$.
Multiplying  (17) by  $\eta_1(1)$  and using (16) we obtain for
the components of $X$ and  $Y$
$$
\eta_1(1)x_j=\eta_1(1)y_j\qquad\qquad(j=1,\ldots,m).\tag18
$$
Therefore, if the $1$-cocycles $\frak{f}$ and $\frak{g}$ are
cohomologous then (18) holds.

Conversely, suppose that (18) holds.  Then it is not difficult to
construct  vectors $Z_2,Z_3,Z_4$ that satisfy (9) and (10), which
is equivalent to (8), i.e. the $1$-cocycles $\frak{f}$ and
$\frak{g}$ are cohomologous. It follows that by going from a
cocycle to an element of the cohomology group, we need  to change
each component in $X$ by $\beta=\alpha\cdot p\m1$ modulo the group
$R$, where $\alpha\in R$. Moreover, if $\eta_1\cdot\beta\in R$,
then must change $\beta$ to $0$. \rightline{\text{\qed}}
\enddemo

\proclaim {Theorem 5} Let $\varepsilon\in R$, where
$\varepsilon^p=1$ and $p>2$. Then the description of the non-split
indecomposable groups $\cry{G,M,\frak{f}}$ is a wild type problem.
\endproclaim

\demo{Proof} For arbitrary matrices $A,B\in M(m,R)$ the map
$$
\Gamma_{A,B}: \quad a\mapsto \left(\smallmatrix
E& 0&E&A&E\\
&\varepsilon E& E&E&B\\
&&\varepsilon^2 E&0&0\\
&&&\varepsilon^3 E&0\\
&&&&\varepsilon^4 E\\
\endsmallmatrix\right)
$$
is an $R$-representation of $G$ of degree $l=5m$. The
$R$-representations $\Gamma_{A,B}$ and $\Gamma_{A_1,B_1}$ are
$R$-equivalent if and only if
$$
C\m1AC\equiv A_1 \pmod{(1-\varepsilon)},\quad \quad \quad
C\m1BC\equiv B_1 \pmod{(1-\varepsilon)}
$$
for some invertible matrix $C$. It follows that the description of
the $R$-representations  $\Gamma_{A,B}$ of $G$ is a wild type
problem.

For the module  affording  the  representation $\Gamma_{A,B}$ of
$G$ we put $R^l$.  Let $X$ be an $m$-dimensional vector over $F$
with $pX\in R^m$. Then there is a $1$-cocycle $\frak{f}_X:
G\longrightarrow \widehat{R}^l$, such that
$\frak{f}_X(a)=(X,0,\ldots,0)+R^l$ . The $1$-cocycles
$\frak{f}_{X}$ and $\frak{f}_Y$ are cohomologous if and only if
$$
(1-\varepsilon)(X-Y)\in R^m.
$$
Putting $X=(p\m1,0,\ldots,0)$ we obtain that
$H^1(G,\widehat{R}^l)\not=0$.\qquad {\text{\qed}}
\enddemo

\enddocument